\documentclass[10pt]{article}
\usepackage{latexsym}
\usepackage{amsfonts}
\title{Orbits of the hyperoctahedral group as Euclidean designs}

\author{B\'ela Bajnok\\ Gettysburg College\\ bbajnok@gettysburg.edu}

\newtheorem{thm}{Theorem}
\newtheorem{defin}[thm]{Definition}
\newtheorem{lem}[thm]{Lemma}
\newtheorem{cor}[thm]{Corollary}
\newtheorem{prop}[thm]{Proposition}

\newcommand{\bc}[2]{{{#1}\choose{#2}}}

\begin{document}

\renewcommand{\today}{September 16, 2006}

\maketitle

\begin{center}

\end{center}

\begin{abstract}

The hyperoctahedral group $H$ in $n$ dimensions (the Weyl group of Lie type $B_n$) is the subgroup of the orthogonal group generated by all transpositions of coordinates and reflections with respect to coordinate hyperplanes.  With ${\bf e_1}, \dots, {\bf e_n}$ denoting the standard basis vectors of $\mathbb{R}^n$ and letting ${\bf x_k}={\bf e_1}+ \cdots + {\bf e_k}$ ($k=1,2,\dots,n$), the set $${\cal I}^n_k={\bf x_k}^H=\{ {\bf x_k}^g \mbox{   } | \mbox{  } g \in H \}$$ is the vertex set of a generalized regular hyperoctahedron in $\mathbb{R}^n$.  

A finite set ${\cal X} \subset \mathbb{R}^n$ with a weight function $w: {\cal X} \rightarrow \mathbb{R}^+$ is called a Euclidean $t$-design, if $$\sum_{r \in R} W_r \overline{f}_{S_{r}} = \sum_{{\bf x} \in {\cal X}} w({\bf x}) f({\bf x})$$ holds for every polynomial $f$ of total degree at most $t$; here $R$ is the set of norms of the points in ${\cal X}$,  $W_r$ is the total weight of all elements of ${\cal X}$ with norm $r$, $S_r$ is the $n$-dimensional sphere of radius $r$ centered at the origin, and $\overline{f}_{S_{r}}$ is the average of $f$ over $S_{r}$.

Here we consider Euclidean designs which are supported by orbits of the hyperoctahedral group.  Namely, we prove that any Euclidean design on a union of generalized hyperoctahedra has strength (maximum $t$ for which it is a Euclidean design) equal to 3, 5, or 7.  We find explicit necessary and sufficient conditions for when this strength is 5 and for when it is 7.  In order to establish our classification, we translate the above definition of Euclidean designs to a single equation for $t=5$, a set of three equations for $t=7$, and a set of seven equations for $t=9$.    

Neumaier and Seidel (1988), as well as Delsarte and Seidel (1989), proved a Fisher-type inequality $|{\cal X}| \geq N(n,p,t)$ for the minimum size of a Euclidean $t$-design in $\mathbb{R}^n$ on $p=|R|$ concentric spheres (assuming that the design is antipodal if $t$ is odd).  A Euclidean design with exactly $N(n,p,t)$ points is called tight.  We exhibit new examples of antipodal tight Euclidean designs, supported by orbits of the hyperoctahedral group, for $N(n,p,t)=$(3,2,5), (3,3,7), and (4,2,7).

\end{abstract}



\section{Introduction}

A Euclidean design is a finite weighted set of points in the $n$-dimensional real Euclidean space $\mathbb{R}^n$ with a certain approximation property, as explained below.  First we introduce a few notations and discuss some background.     

The norm of a point ${\bf x}\in \mathbb{R}^n$, denoted by $|| {\bf x} ||$, is its distance from the origin; the collection of all points with given norm $r>0$ is the sphere $S_r^{n-1}$.  For a finite set ${\cal X} \subset \mathbb{R}^n \setminus \{ {\bf 0} \}$, we call $R=\{ || {\bf x} ||  \mbox{   } | \mbox{   } {\bf x} \in {\cal X} \}$ the \emph{norm spectrum} of ${\cal X}$.  (In this paper, for convenience, we exclude the possibility of ${\bf 0} \in {\cal X}$; see \cite{BanBan:Prepa} for a discussion.)  We can partition ${\cal X}$ into \emph{layers} ${\cal X} = \cup_{r \in R} {\cal X}_{r}$ where ${\cal X}_r={\cal X} \cap S_{r}^{n-1}$.   A weight function on ${\cal X}$ is a function $w: {\cal X} \rightarrow \mathbb{R}^+$; the \emph{weight distribution of $w$ on $R$} is the function $W: R \rightarrow \mathbb{R}^+$ given by $W_r=\sum_{{\bf x} \in {\cal X}_r} w({\bf x})$.   (Throughout this paper, we only consider positive weights.) 

We denote the spaces of real polynomials, homogeneous polynomials, and homogeneous harmonic polynomials on $n$ variables by $\mathrm{Pol}( \mathbb{R}^n)$, $\mathrm{Hom}( \mathbb{R}^n)$, and $\mathrm{Harm}( \mathbb{R}^n)$, respectively.  Often we will restrict the domain of these polynomials to a subset ${\cal Y}$ of $\mathbb{R}^n$ or their degrees to a fixed integer $s$; the corresponding polynomial spaces will be denoted by $\mathrm{Pol}_s( {\cal Y})$, etc.  

Let $\sigma^n$ denote the regular surface measure on $S_r^{n-1}$.  The (continuous) \emph{average value} of a polynomial $f$ on $S_r^{n-1}$ is $$\overline{f}_{S_r^{n-1}}=\frac{1}{\sigma^n(S_r^{n-1})} \int_{S_r^{n-1}} f ({\bf x}) d \sigma^n({\bf x}),$$ where 
$$\sigma^n(S_r^{n-1}) =\frac{2 r^{n-1} \pi^{\frac{n}{2}}}{\Gamma({\frac{n}{2}})}$$ is the surface area of $S_r^{n-1}$.

We are now ready to state our definition of Euclidean $t$-designs.

\begin{defin}[Neumaier--Seidel, \cite{NeuSei:1988a}]\label{DEF}

Let ${\cal X} \subset \mathbb{R}^n \setminus \{ {\bf 0} \}$ be a finite set with norm spectrum $R$, and suppose that a weight function $w$ is given on ${\cal X}$ which has weight distribution $W$ on $R$.  Let $t$ be non-negative integer.  We say that $({\cal X},w)$ is a \textup{\textbf{Euclidean $t$-design}} if 
$$\sum_{r \in R} W_r \overline{f}_{S_{r}^{n-1}} = \sum_{{\bf x} \in {\cal X}} w({\bf x}) f({\bf x})$$ holds for every $f \in \mathrm{Pol}_t( \mathbb{R}^n)$.  The largest value of $t$ for which $({\cal X},w)$ is a Euclidean $t$-design is called the \textup{\textbf{(maximum) strength}} of the design.

\end{defin} 

An important special case of Euclidean designs is spherical designs: there we assume that all points in the design have the same norm, and we do not allow the points to have different weights.  Spherical designs enjoy a vast and rapidly growing literature, and have been studied from a variety of perspectives, including algebra, combinatorics, functional analysis, geometry, number theory, numerical analysis, and statistics.  For general references, please see \cite{Ban:1988a}, \cite{ConSlo:1999a}, \cite{DelPac:Prepa}, \cite{DelGoeSei:1977a}, \cite{EriZin:2001a}, \cite{God:1993a}, \cite{GoeSei:1979a}, \cite{GoeSei:1981a}, \cite{Sei:1996a}, and \cite{SeyZas:1984a}.  The concept of Euclidean designs was introduced by Neumaier and Seidel in 1988 in \cite{NeuSei:1988a} as a generalization of spherical designs, and was subsequently studied by Delsarte and Seidel in \cite{DelSei:1989a}, Seidel in \cite{Sei:1990a} and in \cite{Sei:1994a}, and just recently by Bannai and Bannai in \cite{BanBan:Prepa}.

In \cite{Baj:2005a} we provided a recursive construction for Euclidean $t$-designs in $\mathbb{R}^n$.  Namely, we showed how to use certain Gauss--Jacobi quadrature formulae to ``lift'' a Euclidean $t$-design $({\cal A}^{n-1},w^{n-1})$ in $\mathbb{R}^{n-1}$ to a Euclidean $t$-design $({\cal A}^{n},w^{n})$ in $\mathbb{R}^{n}$.  Our recursion preserved both the norm spectrum $R$ and the weight distribution $W$ on $R$; that is, we had $$R=\{ || {\bf x} ||  \mbox{   } | \mbox{   } {\bf x} \in {\cal A}^{n-1} \}=\{ || {\bf x} ||  \mbox{   } | \mbox{   } {\bf x} \in {\cal A}^{n} \}$$ and $$W_r=\sum_{{\bf x} \in {\cal A}^{n-1}_r} w({\bf x})=\sum_{{\bf x} \in {\cal A}^{n}_r} w({\bf x})$$ for each $r \in R$.  Our recursion yields Euclidean $t$-designs in any dimension; we note that this construction for Euclidean $t$-designs is ``semi-explicit'' in the sense that the co{o}rdinates and the weights of the points in the design are given in terms of the roots of certain Gegenbauer polynomials.  Additionally, in \cite{Baj:2005a} we provided explicit constructions for Euclidean designs in the plane ($n=2$); in particular, we analyzed Euclidean $t$-designs of minimum size for every $t$.  Therefore from now on we assume that $n \geq 3$ unless otherwise noted.

The subject of the present paper is a family of very explicit Euclidean designs; namely, we consider point sets which are unions of orbits of single points under the action of the hyperoctahedral group.  We make this more precise as follows.

Let $O(n)$ be the group of isometries of $S^{n-1}$.  For $1 \leq i \leq n$, let $\rho_i \in O(n)$ be the reflection with respect to the coordinate hyperplane $x_{i}=0$; for $1 \leq i_1 < i_2 \leq n$, let $\sigma_{i_1,i_2}$ be the reflection with respect to the hyperplane $x_{i_1}=x_{i_2}$.  Then the sets $\{ \rho_i \mbox{   } | \mbox{   } 1 \leq i \leq n \}$ and $\{ \sigma_{i_1,i_2} \mbox{   } | \mbox{   } 1 \leq i_1 < i_2 \leq n \}$ generate subgroups $H_{\rho}$ and $H_{\sigma}$ of $O(n)$ which are isomorphic to $(\mathbb{Z}/2\mathbb{Z})^n$ and the symmetric group $S_n$, respectively.  The subgroup $H$ generated by all these reflections (the semidirect product of $H_{\rho}$ and $H_{\sigma}$) is called the \emph{hyperoctahedral group} (or the Weyl group of Lie type $B_n$).  We see that $H$ has order $2^n \cdot n!$.

If $g \in H$, then ${\bf x}^g$ denotes the image of ${\bf x} \in \mathbb{R}^n$ after applying $g$; we also set ${\bf x}^H=\{ {\bf x}^g \mbox{   } | \mbox{  } g \in H \}$.  We make the following definition.

\begin{defin}  Let ${\bf e_1}, \dots, {\bf e_n}$ be the standard basis vectors of $\mathbb{R}^n$, and for $k=1,2,\dots,n$, set ${\bf x_k}={\bf e_1}+ \cdots + {\bf e_k}$.  Let $H$ be the hyperoctahedral group of dimension $n$.  The set $${\cal I}^n_k={\bf x_k}^H=\{ {\bf x_k}^g \mbox{   } | \mbox{  } g \in H \}$$ is the vertex set of a \textup{\textbf{generalized regular hyperoctahedron}} in $\mathbb{R}^n$.  
\end{defin}

Note that the $n$ generalized hyperoctahedra in dimension $n$, together with the origin, form a partition of the $3^n$ elements in ${\cal I}^n =\{-1,0,1\}^n$ according to their norm.  For example, for $n=3$, this partition is the following:  $${\cal I}^3 =\{ {\bf 0} \} \cup \{ \mbox{octahedron} \}  \cup \{ \mbox{cuboctahedron} \} \cup \{ \mbox{cube} \}.$$  Furthermore, we have $$|{\cal I}^n_k|=2^k \cdot \bc{n}{k}.$$

Our goal is to use such generalized hyperoctahedra to construct Euclidean designs.  We consider sets ${\cal X}$ of the form $${\cal X}={\cal X}(J)=\cup_{k \in J} \frac{r_k}{\sqrt{k}} {\cal I}^n_k$$ where $J \subseteq \{1,2,\dots,n\}$ and $r_k>0$ for every $k \in J$.  Note that ${\cal X}$ has norm spectrum $R=\{r_k \mbox{  } | \mbox{  } k \in J \}$.  The weight function $w: {\cal X} \rightarrow \mathbb{R}^+$ on ${\cal X}$ will be constant on each layer of ${\cal X}$ (see \cite{BanBan:Prepa}); let us denote the weight of ${\bf x} \in \frac{r_k}{\sqrt{k}}{\cal I}^n_k$ by $w({\bf x})=w_k$.  Then $({\cal X(J)},w)$ is an \emph{antipodal design}, that is, for every $x \in {\cal X}$, $ -x \in {\cal X}$ and $w(-x)=w(x)$. 

Our goal is to determine necessary and sufficient conditions for each dimension $n$, norm spectrum $R$, strength $t$, and index set $J$, whether a (positive) weight function exist for which $({\cal X(J)},w)$, as defined above, is a Euclidean $t$-design.  We carry this out here for $|R| \leq 3$ and $|J| \leq 3$.  Note that, since our sets are antipodal, the maximum strength of each Euclidean design supported by them is necessarily odd.

It turns out that, in order to state our results precisely, the following function is useful.  For positive integers $k_1, k_2 \leq n$ we define $$G(k_1,k_2)=(n+2-3k_1)(n+2-3k_2)+6(k_1-1)(k_2-1)+2(n-1).$$  We say that $n$ has \emph{property $G$} if there are positive integers $k_1, k_2 \leq n$ for which $G(k_1,k_2)=0$; in this case we write $n \in G$.  

It is an interesting question to characterize integers with property $G$; these values under a 100 are

\noindent 1, 2, 5, 8, 10, 11, 14, 16, 17, 20, 23, 26, 28, 29, 31, 32, 35, 38, 41, 44, 46, 47, 50, 52, 53, 56, 59, 61, 62, 64, 65, 68, 71, 73, 74, 76, 77, 80, 82, 83, 86, 89, 91, 92, 94, 95, 98, 100.  

\noindent Clearly, no $n$ value with $n \equiv 0$ mod 3 has property $G$.  We can also see that all $n \equiv 2$ mod 3 has property $G$: for example, we can take $k_1=1$ and $k_2=\frac{n+4}{3}$.  Among integers congruent to 1 mod 3, there are 176 such integers under a thousand (slightly over half).  We can also find infinite subsequences of $n \equiv 1$ mod 3 with property $G$: for example, each $n$ of the form $n=90c^2+45c+1$ has property $G$ (with $k_1=3c+1$ and $k_2=30c^2+17c+2$).

We now state our characterization as follows.  

\begin{thm}\label{GGA}

Let ${\cal X}(J)=\cup_{k \in J} \frac{r_k}{\sqrt{k}} {\cal I}^n_k$ be a union of generalized hyperoctahedra, as defined above.

\begin{enumerate}

\item The strength of $({\cal X}(J),w)$ is 3, 5, or 7 for all choices of $R$, $J$, and $w$.

\item Suppose that $|R|=1$.

\begin{description}

\item{$|J|=1; t=5$} \\ If $J=\{k\}$, then there is a weight function $w$ for which $({\cal X}(J),w)$ a Euclidean $5$-design, if and only if, $n \equiv 1$ mod 3 and $k=\frac{n+2}{3}$.

\item{$|J|=1; t=7$} \\ If $|J|=1$, then $({\cal X}(J),w)$ is never a Euclidean $7$-design.

\item{$|J|=2; t=5$} \\ Let $J=\{k_1,k_2\}$ where $k_1 < k_2$.  Then there is weight function $w$ for which $({\cal X}(J),w)$ is a Euclidean $5$-design, if and only if, $k_1 < \frac{n+2}{3} < k_2$.

\item{$|J|=2; t=7$} \\ If $J=\{k_1,k_2\}$, then there is a weight function $w$ for which $({\cal X}(J),w)$ is a Euclidean $7$-design, if and only if, $G(k_1,k_2)=0$.

\item{$|J|=3; t=7$} \\ Let $J=\{k_1,k_2,k_3\}$ where $k_1<k_2<k_3$.  Then there is a weight function $w$ for which $({\cal X}(J),w)$ is a Euclidean $7$-design, if and only if, $G(k_1,k_2)>0$, $G(k_2,k_3)>0$, and $G(k_1,k_3)<0$.

\end{description}

\item Suppose that $|R|=2$.

\begin{description}

\item{$|J|=2; t=5$} \\   Let $J=\{k_1,k_2\}$ where $k_1 < k_2$.  Then there is weight function $w$ for which $({\cal X}(J),w)$ is a Euclidean $5$-design, if and only if, $k_1 < \frac{n+2}{3} < k_2$.

\item{$|J|=2; t=7$} \\  If $|J| =2$, then $({\cal X}(J),w)$ is never a Euclidean $7$-design.  

\item{$|J|=3; t=7$} \\  Let $J=\{k_1,k_2,k_3\}$ where $k_1<k_2<k_3$.  Then there is a weight function $w$ for which $({\cal X}(J),w)$ is a Euclidean $7$-design, if and only if, $n \equiv 1$ mod 3, $k_2=\frac{n+2}{3}$, $r_{k_1}=r_{k_3}$, and $G(k_1,k_3)<0$.

\end{description}

\item Suppose that $|R|=3$.

\begin{description}

\item{$|J|=3; t=7$} \\  Let $J=\{k_1,k_2,k_3\}$ where $k_1<k_2<k_3$.  Then there is a weight function $w$ for which $({\cal X}(J),w)$ is a Euclidean $7$-design, if and only if, $G(k_1,k_2)>0$, $G(k_2,k_3)>0$, $G(k_1,k_3)<0$, and \begin{equation}\label{Q} 
\sum_{i=1,2,3} k_i \cdot (n+2-3k_i) \cdot (k_{i+1}-k_{i+2}) \cdot G(k_{i+1},k_{i+2}) \cdot \frac{1}{r_{k_{i}}^2} =0,
\end{equation} holds where the indeces are to be taken mod 3 (i.e. $k_4=k_1$ and $k_5=k_2$).

\end{description}

\end{enumerate}

\end{thm}

We prove Theorem \ref{GGA} in Section 3.  As a consequence, we have the following.

\begin{cor}\label{GGB}

Let $\tau(p,j)$ be the maximum achievable strength of a Euclidean design formed by the union of $j$ generalized hyperoctahedra on $p$ concentric spheres.  We have the following.

$$\begin{array}{lll}

 \tau(1,1)= \left\{ 
\begin{array}{l}
3 \mbox{   if   }  n \not \equiv 1 (3) \\ \\
5 \mbox{   if   }  n \equiv 1 (3) 
\end{array}
\right.

&
\mbox{         }
&

\tau(2,2)= 5 \mbox{   for all   }  n

\\ \\ \\

  \tau(1,2)= \left\{ 
\begin{array}{l}
5 \mbox{   if   }  n \not \in G \\ \\
7 \mbox{   if   }  n \in G 
\end{array}
\right.

&
\mbox{         }
&

 \tau(2,3)= \left\{ 
\begin{array}{l}
5 \mbox{   if   }  n \not \equiv 1 (3) \\ \\
7 \mbox{   if   }  n \equiv 1 (3) 
\end{array}
\right.

\\ \\ \\

 \tau(1,3)= 7 \mbox{   for all   }  n

&
\mbox{         }
&

 \tau(3,3)= \left\{ 
\begin{array}{l}
5 \mbox{   if   }  n =4 \\ \\
7 \mbox{   if   }  n \not =4  
\end{array}
\right.

\end{array}$$

\end{cor}

In particular, we see that

\begin{itemize}
\item for every $n \geq 3$ and for $p \in \{1,2\}$, we can find at most two generalized hyperoctahedra whose union supports a Euclidean $5$-design;  

\item for every $n \geq 3$ other than $n=4$ and for $p \in \{1,3\}$, we can find at most three generalized hyperoctahedra whose union supports a Euclidean $7$-design; and 

\item for $p=2$, we can find three generalized hyperoctahedra whose union supports a Euclidean $7$-design only if $n \equiv 1$ mod 3.

\end{itemize}

In order to establish our results, we introduce a necessary and sufficient criterion for $({\cal X},w)$ to be a Euclidean $t$-design which, in our situation, is much easier to use than Definition \ref{DEF}.  First we make the following definition.

\begin{defin}\label{DEFF}

Let $H$ be the hyperoctahedral group in dimension $n$.  A subset ${\cal X}$ of $\mathbb{R}^n$ is said to be \textup{\textbf{fully symmetric}} if ${\cal X}^g={\cal X}$ for every $g \in H$.  Furthermore, a Euclidean design $({\cal X},w)$ is fully symmetric if ${\cal X}$ is a fully symmetric set, and $w({\bf x})=w({\bf y})$ whenever ${\bf x}={\bf y}^g$ for some $g \in H$. 

\end{defin}
Note that the designs supported on generalized hyperoctahedra are fully symmetric.

The following result proves to be quite useful.

\begin{thm}\label{FJ} Suppose that $({\cal X},w)$ is fully symmetric with norm spectrum $R$, and assume that its weight function is a constant $w_r$ on each layer ${\cal X}_r$ of ${\cal X}$.  Define polynomials $f_{4,2}$, $f_{6,3}$, $f_{8,2}$, and $f_{8,4}$ as follows.
\begin{eqnarray*} 
f_{4,2} & = & x_1^4-6x_1^2x_2^2+x_2^4, \\ \\
f_{6,3} & = & 2(x_1^6+x_2^6+x_3^6)-15(x_1^4x_2^2+x_1^2x_2^4+x_1^4x_3^2+x_1^2x_3^4+x_2^4x_3^2+x_2^2x_3^4)+180x_1^2x_2^2x_3^2,\\ \\
f_{8,2} & = & x_1^8-28x_1^6x_2^2+70x_1^4x_2^4-28x_1^2x_2^6+x_2^8, \\ \\
f_{8,4} & = & 3(x_1^8+x_2^8+x_3^8+x_4^8)-28(x_1^6x_2^2+\cdots+x_3^2x_4^6)+210(x_1^4x_2^2x_3^2+\cdots+x_2^2x_3^2x_4^4) -3780x_1^2x_2^2x_3^2x_4^2.
\end{eqnarray*}  
Then we have the following.
\begin{enumerate}

\item $({\cal X},w)$ is at least a Euclidean $3$-design.

\item $({\cal X},w)$ is a Euclidean $5$-design, if and only if, 
\begin{equation}\label{DA}
\sum_{r \in R} \sum_{{\bf x} \in {\cal X}_r} w_r f_{4,2}({\bf x}) =0.
\end{equation}  

\item $({\cal X},w)$ is a Euclidean $7$-design, if and only if, 
\begin{equation}\label{DB}
 \left\{ 
\begin{array}{l}
\sum_{r \in R} \sum_{{\bf x} \in {\cal X}_r} w_r r^{2s_1} f_{4,2}({\bf x}) =0 \\ \\
\sum_{r \in R} \sum_{{\bf x} \in {\cal X}_r} w_r f_{6,3}({\bf x}) =0 \\ \\ 
\end{array}
\right.
\end{equation} holds for $s_1=0, 1$.

\item $({\cal X},w)$ is a Euclidean $9$-design, if and only if, 
\begin{equation}\label{DC}
 \left\{ 
\begin{array}{l}
\sum_{r \in R} \sum_{{\bf x} \in {\cal X}_r} w_r r^{2s_1} f_{4,2}({\bf x}) =0 \\ \\
\sum_{r \in R} \sum_{{\bf x} \in {\cal X}_r} w_r r^{2s_2} f_{6,3}({\bf x}) =0 \\ \\
\sum_{r \in R} \sum_{{\bf x} \in {\cal X}_r} w_r f_{8,2}({\bf x}) =0 \\ \\
\sum_{r \in R} \sum_{{\bf x} \in {\cal X}_r} w_r f_{8,4}({\bf x}) =0
\end{array}
\right.
\end{equation} holds for $s_1=0, 1, 2$ and $s_2=0, 1$.

\end{enumerate}

\end{thm} 
We prove Theorem \ref{FJ} in the next section.

Let us now examine Euclidean designs of minimum size.  For a non-negative integer $s$, let 
$$d(s):=\dim \mathrm{Hom}_s( \mathbb{R}^n) = \bc{s+n-1}{n-1}.$$
We have the following \emph{Fisher-type inequality}.  The inequality provides the minimum size for the case of an even $t$; for odd $t$ we only have a lower bound if the design is antipodal.

\begin{thm}[Delsarte--Seidel, \cite{DelSei:1989a}]\label{FB}  For a positive integer $k$, let $$N_k=d \left( \left\lfloor \frac{t}{2} \right\rfloor +2-2k \right) + d \left( \left\lfloor \frac{t-1}{2} \right\rfloor +2-2k \right)$$ and $$N(n,p,t)=\sum_{k=1}^{p} N_k.$$  Suppose that $({\cal X},w)$ is a Euclidean $t$-design on $p$ layers in $\mathbb{R}^n$; if $t$ is odd, assume further that the design is antipodal.
Then we have $|{\cal X}| \geq N(n,p,t).$ 

\begin{defin}\label{JJ}

We say that a Euclidean $t$-design on $p$ concentric spheres is \textup{\textbf{tight}} if it has size $N(n,p,t)$.

\end{defin}

\end{thm}

In \cite{Baj:2005a} we explicitly constructed tight Euclidean designs in the plane for every $t$ and every $p \leq \left\lfloor \frac{t+5}{4}\right\rfloor$.  Our construction was a generalization of the example for $p=2$ and $t=4$ which already appeared in \cite{BanBan:Prepa}; namely, we have shown how a union of concentric regular polygons, with appropriate weights and rotations about the origin, forms a tight Euclidean design.    

Let us now turn to the case of tight Euclidean designs in higher dimensions ($n \geq 3$).  Since here we restrict our search to fully symmetric designs which are antipodal, we may assume that $t$ is odd.  

It is easy to see that for $t=1$ a pair of two antipodal points, with equal weights, forms an antipodal tight 1-design in $\mathbb{R}^n$; and that for $t=3$, the $n$ pairs $\pm r_i \cdot {\bf e_i}$, with weights inversely proportional to their squared norms, form antipodal tight 3-designs in $\mathbb{R}^n$ (in particular, the $2n$ vertices of the octahedron ${\cal I}^n_1$ form an antipodal tight Euclidean 3-design).  Bannai \cite{Ban:Prepb} classified all antipodal tight Euclidean 3-designs and all antipodal tight Euclidean 5-designs with $p=2$.  Here we provide examples of antipodal tight Euclidean designs for $(n,p,t)=$(3,2,5), (3,3,7), and (4,2,7).  Namely, we will prove that 

\begin{enumerate}
\item
in $\mathbb{R}^3$, the union of an octahedron and a cube, with appropriate weights, forms an antipodal tight 5-design;
\item
in $\mathbb{R}^3$, the union of an octahedron, a cuboctahedron, and a cube, with appropriate weights, forms an antipodal tight 7-design; and
\item
in $\mathbb{R}^4$, the points of minimum non-zero norm in the lattice $D_4$ together with the points of minimum non-zero norm in the dual lattice $D_4^*$, with appropriate weights, form an antipodal tight 7-design.
\end{enumerate}

It might be useful to summarize these examples more precisely, as follows.
 
\begin{cor}\label{FG}

Let $({\cal X}(J)=\cup_{k \in J} \frac{r_k}{\sqrt{k}} {\cal I}^n_k,w)$ be a Euclidean design supported by a union of generalized hyperoctahedra, as defined above.  Let $r, \rho, w \in \mathbb{R}^+$, $\rho \not = r$.

\begin{enumerate} 

\item For $n=3$, $p=2$, and $t=5$, choose $({\cal X}(J),w)$ according to the following table.
$$\begin{array}{||c||c|c||}
\hline \hline  
k \in J & 1 & 3  \\ \hline \hline
r_k & r & \rho \\ \hline
w_k & w & \frac{9}{8} \frac{r^4}{\rho^4} \cdot w \\  \hline \hline
\end{array}$$ 
Then $({\cal X}(J),w)$ is an antipodal tight Euclidean 5-design of size 14 in $\mathbb{R}^3$.

\item For $n=3$, $p=3$, and $t=7$, choose $({\cal X}(J),w)$ according to the following table.
$$\begin{array}{||c||c|c|c||}
\hline \hline   
k \in J & 1 & 2 & 3  \\ \hline \hline
r_k & r & \sqrt{\frac{3r^2+2\rho^2}{5}} \cdot \frac{r}{\rho} & \sqrt{\frac{3r^2+2\rho^2}{5}} \\ \hline
w_k & w & \frac{100\rho^6}{(3r^2+2\rho^2)^3} \cdot w & \frac{675 r^6}{8(3r^2+2\rho^2)^3} \cdot w \\ \hline \hline 
\end{array}$$ 
Then $({\cal X}(J),w)$ is an antipodal tight Euclidean 7-design of size 26 in $\mathbb{R}^3$.

\item For $n=4$, $p=2$, and $t=7$, choose $({\cal X}(J),w)$ according to the following table.
$$\begin{array}{||c||c|c|c||}
\hline \hline  
k \in J & 1 & 2 & 4  \\  \hline \hline
r_k & r & \rho & r \\ \hline 
w_k & w & \frac{r^6}{\rho^6} \cdot w & w \\  \hline \hline
\end{array}$$ 
Then $({\cal X}(J),w)$ is an antipodal tight Euclidean 7-design of size 48 in $\mathbb{R}^4$.

\end{enumerate}

\end{cor}

We prove Corollary \ref{FG} in Section 4.  In particular, we will see that each design has the same stated strength even when $\rho=r$; however, they are only tight when $\rho \not =r$.  Note also that in (2) above $\rho \not =r$ implies $|\{r_1,r_2,r_3\}| =3$ (while $\rho =r$ would yield $r_1=r_2=r_3$).

Finally, we point out that, like the results in \cite{Baj:2005a}, \cite{BanBan:Prepa}, and \cite{Ban:Prepb}, Corollary \ref{FG} disproves the conjecture of Neumaier and Seidel in \cite{NeuSei:1988a} (see also \cite{DelSei:1989a}) that there are no tight Euclidean $t$-designs with $p \geq 2$ and $t \geq 4$.  It seems to be an interesting problem to classify all other tight Euclidean designs.

\section{Harmonic polynomials over $\mathbb{R}^n$}

In this section we prove Theorem \ref{FJ}.   In order to do this, we first review some information on harmonic polynomials over $\mathbb{R}^n$; then develop some very useful results about a special subspace of $\mathrm{Harm}_s( \mathbb{R}^n)$.

There are several equivalent definitions of Euclidean designs.  For our purposes in this paper, the following will be convenient.  Recall that a polynomial is \emph{harmonic} if it satisfies Laplace's equation $\Delta f=0$.  The set of homogeneous harmonic polynomials of degree $s$ in $\mathbb{R}^n$ forms the vector space $\mathrm{Harm}_{s}(\mathbb{R}^n)$ with $$\dim \mathrm{Harm}_{s}(\mathbb{R}^n)=\bc{n+s-1}{n-1}-\bc{n+s-3}{n-1}.$$

\begin{prop}[Neumaier--Seidel, \cite{NeuSei:1988a}]\label{FH}

The weighted set $({\cal X},w)$ is a Euclidean $t$-design in $\mathbb{R}^n$, if and only if,  
$$\sum_{{\bf x} \in {\cal X}} w({\bf x})||{\bf x}||^{2s_1} f({\bf x})=0$$ for every $0 \leq 2s_1 \leq t$ and $f \in  \mathrm{Harm}_s( \mathbb{R}^n)$ with $1 \leq s \leq t-2s_1$.

\end{prop}

An explicit basis for $\mathrm{Harm}_{s}(\mathbb{R}^n)$ can be found using Gegenbauer polynomials, as follows.  Recall that, for fixed $\alpha > -1$, the \emph{Gegenbauer polynomial} (also known as the ultraspherical polynomial) $P_s^{\alpha}: I=[-1,1] \rightarrow \mathbb{R}$  can be defined as $$P_s^{\alpha}(x)=\frac{(-1)^s}{2^ss!(1-x^2)^{\alpha-1/2}}\frac{d^s (1-x^2)^{\alpha+s-1/2}}{dx^s}.$$  (This is one of the various normalizations used.)
The Gegenbauer polynomial $P_s^{\alpha}$ has degree $s$ and has $s$ distinct roots in the interval $(-1,1)$; it is an even function when $s$ is even and an odd function when $s$ is odd.
 
In order to form a basis for $\mathrm{Harm}_{s}(\mathbb{R}^n)$, choose integers $m_0, m_1, \dots, m_{n-2}$ with $0 \leq m_{n-2} \leq \cdots \leq m_1 \leq m_0=s$.  For $k=0,1,\dots, n-3$, set $$r_k=r_k(x_{k+1},\dots,x_{n}) =\sqrt{\sum_{l=k+1}^n x_l^2},$$ and define
$$g_k=g_k(x_{k+1},\dots,x_{n}) =r_k^{m_k-m_{k+1}} P_{m_k-m_{k+1}}^{m_{k+1}+(n-k-2)/2} \left(\frac{x_{k+1}}{r_k} \right)$$ where $P_{m_k-m_{k+1}}^{m_{k+1}+(n-k-2)/2}$ is the Gegenbauer polynomial defined above.  Note that $g_k$ is a homogeneous polynomial of degree $m_k-m_{k+1}$.

Let also $$h_1(x_{n-1},x_n)=\mathrm{Re}(x_{n-1}+ix_n)^{m_{n-2}}$$ and $$h_2(x_{n-1},x_n)=\mathrm{Im}(x_{n-1}+ix_n)^{m_{n-2}}.$$

Finally, for integer(s) $1 \leq \mu \leq \mathrm{min} \{2, m_{n-2}+1\}$, define 
$$f_{m_0,m_1,\dots,m_{n-2},\mu}(x_1,\dots,x_n)=h_{\mu}(x_{n-1},x_n) \cdot \prod_{k=0}^{n-3}g_k(x_{k+1},\dots,x_{n})$$
and the set $$\Phi_s^n=\{f_{m_0,m_1,\dots,m_{n-2},\mu}(x_1,\dots,x_n) \}.$$
Note that $f_{m_0,m_1,\dots,m_{n-2},\mu}(x_1,\dots,x_n)$ is a homogeneous polynomial of degree $s$ and that 
\begin{eqnarray*}
|\Phi_s^n| & = & |\{  (m_{n-2},\dots, m_0,\mu) \in \mathbb{Z}^n \mbox{   } | \\ & & \mbox{             } \mbox{             } 0 \leq m_{n-2} \leq \cdots \leq m_0=s, 1 \leq \mu \leq \mathrm{min} \{2, m_{n-2}+1\} \}| \\
& = & 2 \cdot \bc{n+s-3}{n-2} + \bc{n+s-3}{n-3} \\
& = & \bc{n+s-1}{n-1}-\bc{n+s-3}{n-1} \\
& = & \dim \mathrm{Harm}_{s}(\mathbb{R}^n).
\end{eqnarray*} Furthermore, we have the following.

\begin{prop}[\cite{Erd:1953a}]\label{AABB}

With the above notations, the set 
$\Phi_s^n$
forms a basis for $\mathrm{Harm}_{s}(\mathbb{R}^n)$.
\end{prop}

For larger values of $n$, Propositions \ref{FH} and \ref{AABB} are not convenient due to the large size of $\Phi_s^n$.  However, if we consider only fully symmetric designs, then the necessary criteria can be greatly reduced as explained below.  

Recall that, as defined by Definition \ref{DEFF}, a Euclidean design $({\cal X},w)$ is fully symmetric if ${\cal X}^g={\cal X}$ for every $g \in H$ and $w({\bf x})=w({\bf y})$ whenever ${\bf x}={\bf y}^g$ for some $g \in H$.  Here $H$ is the hyperoctahedral group in dimension $n$, that is the subgroup of the orthogonal group $O(n)$ generated by all transpositions and reflections with respect to coordinate hyperplanes.  For a polynomial $f \in \mathrm{Pol}( \mathbb{R}^n)$, we let $f^g$ denote the function $f({\bf x}^g)$.  As before, for $1 \leq i \leq n$, we let $\rho_i \in O(n)$ be the reflection with respect to the coordinate hyperplane $x_{i}=0$, and for $1 \leq i_1 < i_2 \leq n$, let $\sigma_{i_1,i_2}$ be the reflection with respect to the hyperplane $x_{i_1}=x_{i_2}$.  Note that $\rho_i$ acts on a polynomial $f$ by replacing its variable $x_i$ be $-x_i$, and $\sigma_{i_1,i_2}$ acts by switching $x_{i_1}$ and $x_{i_2}$.

Let us first address symmetry with respect to coordinate hyperplanes.  For this purpose, we are specifically interested in \emph{fully even} harmonic polynomials.

\begin{defin}

A polynomial $f \in \mathbb{R}[x_1,\dots,x_n]$ is said to be \textup{\textbf{fully even}} if $f^{\rho_i}=f$ for every $1 \leq i \leq n$, and \textup{\textbf{partially odd}} if $f^{\rho_i}=-f$ for some $1 \leq i \leq n$.  We let $\mathrm{FEvenHarm}_{s}(\mathbb{R}^n)$ and $\mathrm{POddHarm}_{s}(\mathbb{R}^n)$ denote the set of fully even polynomials and the set of partially odd polynomials in $\mathrm{Harm}_{s}(\mathbb{R}^n)$, respectively. 

\end{defin}

Note that in a polynomial which belongs to $\mathrm{FEvenHarm}_{s}(\mathbb{R}^n)$, in every term every variable has an even degree; while in members of $\mathrm{POddHarm}_{s}(\mathbb{R}^n)$, at least one variable appears only with an odd degree in every term.

We let $$\mathrm{FEven}\Phi_s^n=\Phi_s^n \cap \mathrm{FEvenHarm}(\mathbb{R}^n),$$ and $$\mathrm{POdd}\Phi_s^n=\Phi_s^n \cap \mathrm{POddHarm}(\mathbb{R}^n).$$   

Recall that a Gegenbauer polynomial of degree $s$ is an even function when $s$ is even and an odd function when $s$ is odd.  Therefore the polynomial $g_k=g_k(x_{k+1},\dots,x_{n})$, defined above, has its variable $x_{k+1}$ with even exponents only if its degree $m_k-m_{k+1}$ is even, and odd exponents only if $m_k-m_{k+1}$ is odd; the variables $x_{k+2}, \dots, x_n$ all appear with even exponents only.  Consequently, we can determine easily which members of $\Phi_s^n$ are in $\mathrm{FEven}\Phi_s^n$, and we can count the elements of $\mathrm{FEven}\Phi_s^n$, as follows.

\begin{prop}\label{HA}  With the notations above, we have the following.
\begin{enumerate}
\item $\mathrm\Phi_s^n = \mathrm{FEven}\Phi_s^n \cup \mathrm{POdd}\Phi_s^n;$

\item
$\mathrm{FEven}\Phi_s^n$ is a basis for $\mathrm{FEvenHarm}_{s}(\mathbb{R}^n)$;

\item A polynomial $f_{m_0,m_1,\dots,m_{n-2},\mu}(x_1,\dots,x_n) \in \Phi_s^n$ is in $\mathrm{FEven}\Phi_s^n$, if and only if, $m_{0},\cdots, m_{n-2}$ are all even and $\mu=1$;
\item If $s$ is even, then
\begin{equation}\label{AAAC}
\dim \mathrm{FEvenHarm}_{s}(\mathbb{R}^n)=\bc{n+\frac{s}{2}-2}{n-2}.
\end{equation} 
\end{enumerate}
\end{prop}

We can now use Propositions \ref{FH} and \ref{HA} to restate Proposition \ref{FH} for fully symmetric designs as follows.  

\begin{prop}\label{ABBB} A fully symmetric weighted set $({\cal X},w)$ is a Euclidean $t$-design in $\mathbb{R}^n$, if and only if,  
$$\sum_{{\bf x} \in {\cal X}} w({\bf x})||{\bf x}||^{2s_1} f({\bf x})=0$$   holds for every $0 \leq 2s_1 \leq t$ and $f \in  \mathrm{FEvenHarm}_s(\mathbb{R}^n)$ with $s$ even and $2 \leq s \leq t-2s_1$.  
\end{prop}

Note that $\dim \mathrm{FEvenHarm}_{s}(\mathbb{R}^n)$ is substantially smaller than $\dim \mathrm{Harm}_{s}(\mathbb{R}^n)$.  However, rather than working with $\mathrm{FEven}\Phi_s^n$ resulting from the algorithm of Proposition \ref{AABB}, we now attempt to choose a more convenient basis for $\mathrm{FEvenHarm}_s(\mathbb{R}^n)$ which further reduces the number of equations to be checked.  

For $s=2$, by (\ref{AAAC}) we have $\dim \mathrm{FEvenHarm}_{2}(\mathbb{R}^n)=n-1$; we find, for example, that $$F_2=\{x_i^2-x_n^2 \mbox{  } | \mbox{  } 1 \leq i \leq n-1 \}$$ forms a basis for $\mathrm{FEvenHarm}_{2}(\mathbb{R}^n)$.  From this, we immediately see that, if the weighted set $(X, w)$ is fully symmetric, then the equation in Proposition \ref{FH} holds for each $f \in F_2$, and therefore we have Theorem \ref{FJ}, 1: every fully symmetric weighted set $(X, w)$ is a Euclidean 3-design.

Before we turn to finding bases for $s \geq 4$, we introduce a notation.  For a fixed positive integer $j \leq n$, let $S_{n,j}$ denote the set of strictly increasing functions from $\{1,2,\dots,j\}$ to $\{1,2,\dots,n\}$; in other words, permutations of the form
$$g=\left( \begin{array}{cccc}
1 & 2 & \dots & j \\
i_1 & i_2 & \dots & i_j
\end{array} \right)$$  
where $1 \leq i_1 < i_2 < \cdots < i_j \leq n$.  If $f \in \mathbb{R}[x_1,x_2,\dots,x_j]$ is a polynomial on $j$ variables and $g \in S_{n,j}$ is as above, then $f^g$ denotes the polynomial $f(x_{i_1},x_{i_2}, \dots, x_{i_j}) \in \mathbb{R}[x_1,x_2,\dots,x_n]$.  We also write $f^{S_{n,j}}$ for the set $\{f^g \mbox{  } | \mbox{  } g \in S_{n,j}\}$; note that this set has size $\bc{n}{j}$.  

In this paper we focus on $t \leq 9$, thus we must find bases for $\mathrm{FEvenHarm}_s(\mathbb{R}^n)$ for $s=4$, $s=6$, and $s=8$.  According to (\ref{AAAC}), we have 
\begin{eqnarray*} 
\dim \mathrm{FEvenHarm}_{4}(\mathbb{R}^n) & = & \bc{n}{2}, \\ 
\dim \mathrm{FEvenHarm}_{6}(\mathbb{R}^n) & = & \bc{n+1}{3}=\bc{n}{2}+\bc{n}{3}, \mbox{   and}\\ 
\dim \mathrm{FEvenHarm}_{8}(\mathbb{R}^n) & = & \bc{n+2}{4}=\bc{n}{2}+2\bc{n}{3}+\bc{n}{4}.
\end{eqnarray*}
This suggests that we may be able to find polynomials $f_{4,2}, f_{6,2}, f_{8,2} \in \mathbb{R}[x_1,x_2]$, $f_{6,3}, f_{8,3,1}, f_{8,3,2}\in \mathbb{R}[x_1,x_2,x_3]$, and $f_{8,4} \in \mathbb{R}[x_1,x_2,x_3,x_4]$ for which 
\begin{eqnarray*} 
F_4 & = & f_{4,2}^{S_{n,2}}, \\ 
F_6 & = & f_{6,2}^{S_{n,2}} \cup f_{6,3}^{S_{n,3}}, \mbox{   and} \\ 
F_8 & = & f_{8,2}^{S_{n,2}} \cup f_{8,3,1}^{S_{n,3}} \cup f_{8,3,2}^{S_{n,3}} \cup f_{8,4}^{S_{n,4}}
\end{eqnarray*} are bases for $\mathrm{FEvenHarm}_s(\mathbb{R}^n)$ for $s=4$, $s=6$, and $s=8$, respectively.

Indeed, with
\begin{eqnarray*} 
f_{4,2} & = & x_1^4-6x_1^2x_2^2+x_2^4, \\ \\
f_{6,2} & = & x_1^6-15x_1^4x_2^2+15x_1^2x_2^4-x_2^6, \\ \\
f_{6,3} & = & 2(x_1^6+x_2^6+x_3^6)-15(x_1^4x_2^2+x_1^2x_2^4+x_1^4x_3^2+x_1^2x_3^4+x_2^4x_3^2+x_2^2x_3^4)+180x_1^2x_2^2x_3^2,\\ \\
f_{8,2} & = & x_1^8-28x_1^6x_2^2+70x_1^4x_2^4-28x_1^2x_2^6+x_2^8, \\ \\
f_{8,3,1} & = & x_1^8-x_2^8+14(x_1^2x_2^6+x_2^6x_3^2)-14(x_1^6x_2^2+x_1^6x_3^2)+210x_1^4x_2^2x_3^2-210x_1^2x_2^4x_3^2, \\ \\
f_{8,3,2} & = & x_1^8-x_3^8+14(x_1^2x_3^6+x_2^2x_3^6)-14(x_1^6x_3^2+x_1^6x_2^2)+210x_1^4x_2^2x_3^2-210x_1^2x_2^2x_3^4, \\ \\
f_{8,4} & = & 3(x_1^8+x_2^8+x_3^8+x_4^8)-28(x_1^6x_2^2+\cdots+x_3^2x_4^6)+210(x_1^4x_2^2x_3^2+\cdots+x_2^2x_3^2x_4^4) \\ & & -3780x_1^2x_2^2x_3^2x_4^2,
\end{eqnarray*}
we see that the sets $F_4$, $F_6$, and $F_8$ defined above are linearly independent; and, since they have the right cardinality, they form bases for $\mathrm{FEvenHarm}_s(\mathbb{R}^n)$ for $s=4$, $s=6$, and $s=8$, respectively.

Note that if $({\cal X},w)$ is fully symmetric, then the equation of Proposition \ref{ABBB} always holds for \emph{skew-symmetric} polynomials, that is those polynomials $f$ for which $f^{\sigma_{i_1,i_2}}=-f$ for some transposition $\sigma_{i_1,i_2}$; this is the case for $f_{6,2}$, $f_{8,3,1}$, and $f_{8,3,2}$.  This completes the proof of Theorem \ref{FJ}.

\emph{Remarks.}  While we do not attempt to generalize Theorem \ref{FJ} for higher values of $s$, we can indicate what one can expect.  Consider first the following combinatorial identity, often referred to as Vandermonde's Convolution,
\begin{equation}\label{CB}
\bc{n+\frac{s}{2}-2}{n-2}=\sum_{j=2}^{\frac{s}{2}}\bc{\frac{s}{2}-2}{j-2}\bc{n}{j},
\end{equation}
which holds for every even $s$ with $s \geq 4$.  (In this paper we think of $n$ as large, in particular, larger than $s/2$ and thus carry the summation to $s/2$; although (\ref{CB}) is of course valid without this assumption.)  

This identity suggests that for each even $s \geq 4$, one might find polynomials $f_{s,j,i}$, where $2 \leq j \leq \frac{s}{2}$ and $1 \leq i \leq \bc{\frac{s}{2}-2}{j-2}$, so that $$F_s=\bigcup_{i,j} f_{s,j,i}^{S_{n,j}}$$ is a basis for $\mathrm{FEvenHarm}_s(\mathbb{R}^n)$.  This implies that the system of equations in Proposition \ref{ABBB} needs to be checked only for less than $2^{\lfloor t/2 \rfloor}$ equations, a quantity independent of $n$ and substantially smaller than $\dim \mathrm{FEvenHarm}_s(\mathbb{R}^n)$.  We intend to pursue this interesting question in a future study.

\section{Designs supported by generalized hyperoctahedra}

In this section we prove Theorem \ref{GGA}.

First, let us recall our notations.  Let $H$ be the hyperoctahedral group in $n$ dimensions.  Let ${\bf e_1}, \dots, {\bf e_n}$ be the standard basis vectors of $\mathbb{R}^n$, and for $k=1,2,\dots,n$, set ${\bf x_k}={\bf e_1}+ \cdots + {\bf e_k}$.  The set $${\cal I}^n_k={\bf x_k}^H=\{ {\bf x_k}^g \mbox{   } | \mbox{  } g \in H \}$$ is the vertex set of a generalized regular hyperoctahedron in $\mathbb{R}^n$.  We consider sets ${\cal X}$ of the form $${\cal X}={\cal X}(J)=\cup_{k \in J} \frac{r_k}{\sqrt{k}} {\cal I}^n_k$$ where $J \subset \{1,2,\dots,n\}$ and $r_k>0$ for every $k \in J$.  We see that ${\cal X}$ is antipodal and has norm spectrum $R=\{r_k \mbox{  } | \mbox{  } k \in J \}$.  We choose a weight function $w: {\cal X} \rightarrow \mathbb{R}^+$ on ${\cal X}$ which is constant on each layer of ${\cal X}$; let us denote the weight of ${\bf x} \in \frac{r_k}{\sqrt{k}}{\cal I}^n_k$ by $w({\bf x})=w_k$.  

Our goal is to determine for each strength $t$ whether an index set $J$, a norm spectrum $R$, and a weight function $w$ exist for which $({\cal X},w)$ is a Euclidean $t$-design.  Note that $({\cal X(J)},w)$ is fully symmetric, and therefore its strength $t$ is an odd positive integer.  We first prove the following.

\begin{prop}\label{BBBB}  Let $({\cal X(J)},w)$ be the union of weighted generalized hyperoctahedra, as defined above.  Then

\begin{enumerate}

\item $({\cal X(J)},w)$ is at least a Euclidean $3$-design for any $J$, $R$, and $w$.

\item There are no choices of $J$, $R$, and $w$ for which $({\cal X(J)},w)$ is a Euclidean $9$-design.

\end{enumerate}

\end{prop}

{\it Proof.}  The first statement follows immediately from Theorem \ref{FJ}, 1.  For the second statement, we will use (\ref{DC}) of Theorem \ref{FJ}.  In particular, the equation involving $f_{8,2}$ implies that, if $({\cal X(J)},w)$ is a Euclidean 9-design, then 
$$\sum_{k \in J} \sum_{{\bf x} \in {\cal I}_k^n} w_r \cdot \frac{r_k^8}{k^4} \cdot f_{8,2}({\bf x}) =0.$$  However, we find that 
$$\sum_{{\bf x} \in {\cal I}^n_k} f_{8,2}({\bf x})= 2^k \cdot \left[ 2\cdot \bc{n-1}{k-1} - 28 \cdot 2 \cdot \bc{n-2}{k-2} +70 \cdot 2 \cdot \bc{n-2}{k-2} \right],$$ a value which is positive for every $k$, and this is a contradiction.

According to Proposition \ref{BBBB}, the maximum strength of $({\cal X(J)},w)$ is 3, 5, or 7.  Next we determine cases when the strength is 5.

\begin{prop}\label{L}

Let ${\cal X}={\cal X}(J)=\cup_{k \in J} \frac{r_k}{\sqrt{k}} {\cal I}^n_k$ be a union of generalized hyperoctahedra, as defined above.

\begin{enumerate}

\item Suppose that $|J|=1$.  If $n \not \equiv 1$ mod 3, then ${\cal X}(J)$ is never a Euclidean $5$-design.  Suppose that $n \equiv 1$ mod 3 and $J=\{k\}$; then $({\cal X}(J),w)$ is a Euclidean $5$-design, if and only if, $k=\frac{n+2}{3}$. \\ In particular, if $n \equiv 1$ mod 3 and $k=\frac{n+2}{3}$, then $({\cal X}(\{k\}),w)$ is a Euclidean $5$-design for any positive $r_k$ and $w_k$. 

\item  Suppose that $|J| \geq 2$.  Then there is weight function $w$ for which $({\cal X}(J),w)$ is a Euclidean $5$-design, if and only if, there exist $k_1 \in J$ and $k_2 \in J$ for which $k_1 < \frac{n+2}{3} < k_2$.  \\ In particular, $({\cal X}(\{1,n\}),w)$ is a Euclidean $5$-design for any positive $r_1$, $r_n$, $w_1$, and $w_n$ which satisfy $\frac{w_1}{w_n} \cdot \frac{r_1^4}{r_n^4} = \frac{n^2}{2^n}$.

\end{enumerate}

\end{prop}

\emph{Proof.}  To see whether a union of generalized hyperoctahedra is a Euclidean $5$-design, we use Theorem \ref{FJ}, 2.  We compute that $$\sum_{{\bf x} \in {\cal I}^n_k} f_{4,2}({\bf x})= 2\cdot 2^k \cdot \bc{n-1}{k-1} - 6 \cdot 2^k \cdot \bc{n-2}{k-2},$$ and so $({\cal X}(J),w)$ is a 5-design, if and only if, 
\begin{equation}\label{M5}
\sum_{k \in J} w_k \cdot \frac{r_k^4}{k^2} \cdot 2^{k+1} \cdot \bc{n-1}{k-1} \cdot \left[1-3 \cdot \frac{k-1}{n-1} \right]  =0.
\end{equation}  From this our Proposition easily follows.

Now we determine when our construction yields Euclidean $t$-designs for $t=7$.  
Recall that for positive integers $k_1, k_2 \leq n$ we have defined the function $$G(k_1,k_2)=(n+2-3k_1)(n+2-3k_2)+6(k_1-1)(k_2-1)+2(n-1).$$

\begin{prop}\label{HHA}

Let ${\cal X}={\cal X}(J)=\cup_{k \in J} \frac{r_k}{\sqrt{k}} {\cal I}^n_k$ be a union of generalized hyperoctahedra, as defined above.

\begin{enumerate}

\item If $|J|=1$, then ${\cal X}(J)$ is never a Euclidean $7$-design. 

\item  Suppose that $|J| =2$ with $J=\{k_1,k_2\}$.  Then there is a weight function $w$ for which $({\cal X}(J),w)$ is a Euclidean $7$-design, if and only if, $r_{k_1} =r_{k_2}$ and $G(k_1,k_2)=0$.

\item  Suppose that $|J| =3$ with $J=\{k_1,k_2,k_3\}$; assume that $k_1<k_2<k_3$.  
\begin{enumerate}
\item If $|R|=1$, then there is a weight function $w$ for which $({\cal X}(J),w)$ is a Euclidean $7$-design, if and only if, $G(k_1,k_2)>0$, $G(k_2,k_3)>0$, and $G(k_1,k_3)<0$.
\item If $|R|=2$, then there is a weight function $w$ for which $({\cal X}(J),w)$ is a Euclidean $7$-design, if and only if, $n \equiv 1$ mod 3, $k_2=\frac{n+2}{3}$, $r_{k_1}=r_{k_3}$, and $G(k_1,k_3)<0$.
\item If $|R|=3$, then there is a weight function $w$ for which $({\cal X}(J),w)$ is a Euclidean $7$-design, if and only if, $G(k_1,k_2)>0$, $G(k_2,k_3)>0$, $G(k_1,k_3)<0$, and \begin{equation}\label{Q} 
\sum_{i=1,2,3} k_i \cdot (n+2-3k_i) \cdot (k_{i+1}-k_{i+2}) \cdot G(k_{i+1},k_{i+2}) \cdot \frac{1}{r_{k_{i}}^2} =0,
\end{equation} holds where the indeces are to be taken mod 3 (i.e. $k_4=k_1$ and $k_5=k_2$).

\end{enumerate}

\end{enumerate}

\end{prop}

To prove Proposition \ref{HHA}, we will need the following Lemma whose proof is straight forward and will be omitted.

\begin{lem}\label{LEM}
For a positive integer $k$ we introduce the following notations:
$$p_k=k \cdot \left(1-3 \cdot \frac{k-1}{n-1} \right) \mbox{   and   }
q_k=3 \cdot \left(1-15 \cdot \frac{k-1}{n-1} +30 \cdot \frac{(k-1)(k-2)}{(n-1)(n-2)}\right).$$  We then have the following.

\begin{enumerate}

\item If $1 \leq k_1 < k_2 \leq n$, then $$p_{k_1}q_{k_2}-p_{k_2}q_{k_1}=\frac{3(n+8)}{(n-1)^2(n-2)} \cdot (k_1-k_2) \cdot G(k_1,k_2)$$

\item If $1 \leq k_1 < k_2 < k_3 \leq n$, then 
\begin{equation}\label{QA} 
\sum_{i=1,2,3} k_i \cdot (n+2-3k_i) \cdot (k_{i+1}-k_{i+2}) \cdot G(k_{i+1},k_{i+2}) =0,
\end{equation} where the indeces are to be taken mod 3 (i.e. $k_4=k_1$ and $k_5=k_2$).

\item If $1 \leq k_1 < k_2 < k_3 < k_4 \leq n$ and $G(k_2,k_3) \leq 0$, then we have 
\begin{enumerate}

\item $k_2< \frac{n+2}{3}$ and $k_3> \frac{n+2}{3}$

\item $p_2>0$ and $p_3<0$; and

\item $G(k_1,k_2) > 0$ and $G(k_3,k_4) > 0$.

\end{enumerate}
\end{enumerate}

\end{lem}

\emph{Proof of Proposition \ref{HHA}.}  Referring to Theorem \ref{FJ}, 3, we compute that $$\sum_{{\bf x} \in {\cal I}^n_k} f_{6,3}({\bf x})= 6\cdot 2^k \cdot \bc{n-1}{k-1} - 90 \cdot 2^k \cdot \bc{n-2}{k-2}+180 \cdot 2^k \cdot \bc{n-3}{k-3},$$ and so $({\cal X}(J),w)$ is a 7-design, if and only if, 

\begin{equation}\label{DB}
 \left\{ 
\begin{array}{l}
\sum_{k \in J} w_k \cdot \frac{r_k^4}{k^2} \cdot 2^{k+1} \cdot \bc{n-1}{k-1} \cdot \left[1-3 \cdot \frac{k-1}{n-1} \right]  =0 \\ \\
\sum_{k \in J} w_k \cdot \frac{r_k^6}{k^3} \cdot 2^{k+1} \cdot \bc{n-1}{k-1} \cdot 3 \cdot \left[1-15 \cdot \frac{k-1}{n-1} +30 \cdot \frac{(k-1)(k-2)}{(n-1)(n-2)}\right]  =0 \\ \\
\sum_{k \in J} w_k \cdot \frac{r_k^6}{k^2} \cdot 2^{k+1} \cdot \bc{n-1}{k-1} \cdot \left[1-3 \cdot \frac{k-1}{n-1} \right]  =0 
\end{array}
\right.
\end{equation}

With the notations introduced in Lemma \ref{LEM} and with
$$u_k=w_k \cdot \frac{2^{k+1}}{k^3} \cdot \bc{n-1}{k-1},$$
(\ref{DB}) becomes
\begin{equation}\label{DBDB}
 \left\{ 
\begin{array}{l}
\sum_{k \in J} u_k \cdot r_k^4 \cdot p_k =0 \\ \\
\sum_{k \in J} u_k \cdot r_k^6 \cdot q_k =0 \\ \\
\sum_{k \in J} u_k \cdot r_k^6 \cdot p_k =0 
\end{array}
\right.
\end{equation} 

We analyze (\ref{DBDB}), as follows.  First, we see that ${\cal X}(J)$ is never a Euclidean $7$-design for $|J|=1$, since there is no value of $k$ for which $p_k=q_k=0$.  

Suppose now that $|J|=2$ with $J=\{k_1,k_2\}$.  We first note that in this case we must have $p_{k_1} \not =0$ and $p_{k_2} \not =0$, and therefore, if (\ref{DBDB}) holds, then $r_{k_1} =r_{k_2}$, in which case (\ref{DBDB}) is equivalent to 
\begin{equation}\label{DBDBDB}
 \left\{ 
\begin{array}{l}
 u_{k_1} \cdot p_{k_1} + u_{k_2} \cdot p_{k_2}=0 \\ \\
 u_{k_1}  \cdot q_{k_1} + u_{k_2} \cdot q_{k_2}=0 
\end{array}
\right.
\end{equation}
This system has non-zero solutions exactly when $p_{k_1}\cdot q_{k_2} - p_{k_2}\cdot q_{k_1} =0$.  By Lemma \ref{LEM}, 1, this is equivalent to $G(k_1,k_2)=0$; by Lemma \ref{LEM}, 3 (b), we also see that in this case $p_{k_1}$ and $p_{k_2}$ have opposite signs, and therefore we may have $u_{k_1}>0$ and  $u_{k_2}>0$.

Suppose next that $|J|=3$ with $J=\{k_1,k_2,k_3\}$; we also assume that $k_1 < k_2 < k_3$.  

We first show that, if (\ref{DBDB}) holds with non-zero $u_{k_1}$, $u_{k_2}$, and $u_{k_3}$, then $G(k_i,k_j) \not =0$ for $1 \leq i<j \leq 3$.  Indeed, from the last two equations of (\ref{DBDB}) we get   
\begin{equation}\label{DBDBDBDB}
 \left\{ 
\begin{array}{l}
 u_{k_1} \cdot r_{k_1}^6 \cdot (p_{k_3}q_{k_1}-p_{k_1}q_{k_3}) + u_{k_2} \cdot r_{k_2}^6 \cdot (p_{k_3}q_{k_2}-p_{k_2}q_{k_3})=0 \\ \\
u_{k_1} \cdot r_{k_1}^6 \cdot (p_{k_2}q_{k_1}-p_{k_1}q_{k_2}) + u_{k_3} \cdot r_{k_3}^6 \cdot (p_{k_2}q_{k_3}-p_{k_3}q_{k_2})=0 \\ \\
u_{k_2} \cdot r_{k_2}^6 \cdot (p_{k_1}q_{k_2}-p_{k_2}q_{k_1}) + u_{k_3} \cdot r_{k_3}^6 \cdot (p_{k_1}q_{k_3}-p_{k_3}q_{k_1})=0 
\end{array}
\right.
\end{equation}
Therefore, if any $p_{k_i}q_{k_j}-p_{k_j}q_{k_i}$, $1 \leq i<j \leq 3$, equals zero, then they all equal zero; by Lemma \ref{LEM}, 1, this is equivalent to $G(k_1,k_2)=G(k_1,k_3)=G(k_2,k_3)=0$, which would contradict Lemma \ref{LEM}, 3 (c).

We thus see that (\ref{DBDB}) is equivalent to
\begin{equation}\label{DBDBDBDB}
 \left\{ 
\begin{array}{l}
 u_{k_2}=\frac{p_{k_1}q_{k_3}-p_{k_3}q_{k_1}}{p_{k_3}q_{k_2}-p_{k_2}q_{k_3}} \cdot \frac{r_{k_1}^6}{r_{k_2}^6} \cdot u_{k_1}  \\ \\
u_{k_3}=\frac{p_{k_2}q_{k_1}-p_{k_1}q_{k_2}}{p_{k_3}q_{k_2}-p_{k_2}q_{k_3}} \cdot \frac{r_{k_1}^6}{r_{k_3}^6} \cdot u_{k_1} \\ \\
\sum_{i=1,2,3} p_{k_i} \cdot (p_{k_{i+1}}q_{k_{i+2}}-p_{k_{i+2}}q_{k_{i+1}}) \cdot \frac{1}{r_{k_{i}}^2} =0
\end{array}
\right.
\end{equation}
Using Lemma \ref{LEM}, 2, the last equation is equivalent to (\ref{Q}), and the first two can be written as 
\begin{equation}\label{DBDBDBDB}
 \left\{ 
\begin{array}{l}
 u_{k_2}=\frac{k_1-k_3}{k_3-k_2} \cdot \frac{G(k_1,k_3)}{G(k_2,k_3)} \cdot \frac{r_{k_1}^6}{r_{k_2}^6} \cdot u_{k_1} \\ \\
u_{k_3}=\frac{k_2-k_1}{k_3-k_2} \cdot \frac{G(k_1,k_2)}{G(k_2,k_3)} \cdot \frac{r_{k_1}^6}{r_{k_3}^6} \cdot u_{k_1} 
\end{array}
\right.
\end{equation}

Our equations show that $u_{k_1}$, $u_{k_2}$, and $u_{k_3}$ are all positive, if and only if, one of the following two possibilities holds: (i) $G(k_1,k_2)>0$, $G(k_1,k_3)<0$, and $G(k_2,k_3)>0$, or (ii) $G(k_1,k_2)<0$, $G(k_1,k_3)>0$, and $G(k_2,k_3)<0$.  However, by Lemma \ref{LEM}, 3 (c), the second possibility cannot occur.

It remains to see what happens when $|R|<3$.  First note that, by Lemma \ref{LEM}, 2, (\ref{Q}) always holds when $r_{k_1}=r_{k_2}=r_{k_3}$.  If $|R|=2$, then we may assume that for some $i \in \{1,2,3\}$, we have $r_{k_i}=r_{k_{i+1}} \not =r_{k_{i+2}}$ (again, indices are to be taken mod 3).  Comparing (\ref{Q}) and (\ref{QA}), we conclude that $k_{i+2}=\frac{n+2}{3}$.  Since $G(k_1,k_3)<0$, by Lemma \ref{LEM}, 3 (a), we see that this can only happen if $i=3$.  Finally, note that $k_2=\frac{n+2}{3}$ and $k_1<k_2<k_3$ imply, by Lemma \ref{LEM}, 3 (c), that $G(k_1,k_2)>0$ and $G(k_2,k_3)>0$.   

This completes our proof of Proposition \ref{HHA} and hence Theorem \ref{GGA}.

\section{Tight Euclidean designs and other examples}

In this section we examine some examples and prove Corollaries \ref{GGB} and \ref{FG}.  First we prove the following.

\begin{prop}\label{KL}

Let ${\cal X}(J)=\cup_{k \in J} \frac{r_k}{\sqrt{k}} {\cal I}^n_k$ be defined as before.  

\begin{enumerate}

\item

For each of the following sets $J$ and integers $p$, there exists a set $R$ with $|R|=p$ and a weight function $w$ for which $({\cal X}(J),w)$ is a Euclidean $5$-design. 

\begin{enumerate}

\item $p=1$ and $J=\{\frac{n+2}{3} \}$, when $n \equiv 1$ mod 3;

\item $p \in \{1,2\}$ and $J=\{1,n\}$.

\end{enumerate}

\item

For each of the following sets $J$ and integers $p$, there exists a set $R$ with $|R|=p$ and a weight function $w$ for which $({\cal X}(J),w)$ is a Euclidean $7$-design.

\begin{enumerate}

\item $p=1$ and $J=\{k_1,k_2 \}$, when $G(k_1,k_2)=0$; for example, $J=\{1, \frac{n+4}{3} \}$, when $n \equiv 2$ mod 3;

\item $p =1$ and $J=\{1,k,n\}$, with $\frac{n+4}{6} < k < \frac{n+4}{3}$;

\item $p=2$ and $J=\{1, \frac{n+2}{3}, n \}$, when $n \equiv 1$ mod 3;

\item $p =3$ and $J=\{1,k,n\}$, with $\frac{n+4}{6} < k < \frac{n+2}{3}$ or with $k = \frac{n+3}{3}$.

\end{enumerate}

\end{enumerate}
\end{prop}

\emph{Proof.}  Parts 1 (a), 1 (b), and 2 (a) are obvious from Theorem \ref{GGA}.  Next, note that $G(1,n)=-2(n-1)(n-2)$, $G(1,k)=(n-1)(n+4-3k)$, and $G(k,n)=-2(n-1)(n+4-6k)$; so 2 (a) follows from Theorem \ref{GGA}, 2 (e).  Part 2 (b) follows similarly from Theorem \ref{GGA}, 3 (c).  Finally, for part 2 (c), note that with $k_1=1$, $k_2=k$, and $k_3=n$, since $G(1,k)>0$, $G(k,n)>0$, and $G(1,n)<0$, in equation (\ref{Q}), the coefficient of $1/r_{k_{1}}^2$ is negative, the coefficient of $1/{r_{k_{3}}^2}$ is positive, and the coefficient of $1/{r_{k_{2}}^2}$ is positive or negative, depending on $k_2 > \frac{n+2}{3}$ or $k_2 < \frac{n+2}{3}$.  In the first case, an appropriate $r_1$ exists for each choice of $r_k$ and $r_n$; in the second case (which can only happen if $k = \frac{n+3}{3}$), an appropriate $r_n$ exists for each choice of $r_1$ and $r_k$.  This completes our proof.

\emph{Proof of Corollary \ref{GGB}.}  First note that the open interval $(\frac{n+4}{6} , \frac{n+4}{3})$ contains an integer for every $n \geq 3$; and the set $(\frac{n+4}{6} , \frac{n+2}{3}) \cup \{\frac{n+3}{3}\}$ contains an integer for $n=3$ and every $n \geq 5$.  Therefore, all cases of Corollary \ref{GGB} follow from Proposition \ref{KL}, except for $p=3$ and $|J|=3$ when $n=4$ -- a case which can easily be checked. 

Our results can be conveniently summarized in the following table.

\begin{small}
$$\begin{array}{||c||c|c|c|c|c||}
\hline \hline 

& n \equiv 0 (3) & n=4 & n \equiv 1 (3), n \not \in G, n \not = 4 & n \equiv 1 (3), n  \in G & n \equiv 2 (3) \\ \hline
\mbox{First few such   } n & 3,6,9,12,\dots & 4 & 7, 13, 19, 22, 25, 34, \dots & 10, 16, 28, 31, 46, \dots & 5, 8, 11, 14,  \dots \\ \hline \hline

\begin{array}{c|c}
p=1 & \begin{array}{c} |J|=1 \\  |J|=2 \\  |J|=3 \\ \end{array}

\end{array}

&

\begin{array}{c} 3 \\  5 \\  7 \\ \end{array}

&

\begin{array}{c} 5 \\  5 \\  7 \\ \end{array}

&

\begin{array}{c} 5 \\  5 \\  7 \\ \end{array}

&

\begin{array}{c} 5 \\  7 \\  7 \\ \end{array}

&

\begin{array}{c} 3\\  7 \\  7 \\ \end{array}

\\ \hline

\begin{array}{c|c}
p=2 & \begin{array}{c}   |J|=2 \\  |J|=3 \\ \end{array}

\end{array}

&

\begin{array}{c}   5 \\  5 \\ \end{array}

&

\begin{array}{c}   5 \\  7 \\ \end{array}

&

\begin{array}{c}   5 \\  7 \\ \end{array}

&

\begin{array}{c}   5 \\  7 \\ \end{array}

&

\begin{array}{c}   5 \\  5 \\ \end{array}

\\ \hline

\begin{array}{c|c}
p=3 & \begin{array}{c}     |J|=3 \\ \end{array}

\end{array}

&

\begin{array}{c}     7 \\ \end{array}

&

\begin{array}{c}    5 \\ \end{array}

&

\begin{array}{c}     7 \\ \end{array}

&

\begin{array}{c}    7 \\ \end{array}

&

\begin{array}{c}     7 \\ \end{array}

\\ \hline

\end{array}$$  
             
\end{small}

It is interesting to further analyze our results for $n=3$ and for $n=4$.  The following two tables show the maximum strength $\tau$ achievable for a Euclidean design of type $({\cal X(J)},w)$ in 3 and in 4 dimensions on $p$ concentric spheres; we also list the sizes of the designs (we omit braces when indicating the set $J$).  All entries follow from previous results, except for those in the last column for $n=4$; there we find Euclidean 7-designs for $p=2$, if and only if, $r_1=r_3=r_4 \not = r_2$ and for $p=3$, if and only if, $r_3 < r_1=r_2<r_4$ or $r_4 < r_1=r_2<r_3$. 

$$\begin{array}{||c||c|c|c|c|c|c|c||}
\hline \hline 
J & 1 & 2 & 3 & 1,2 & 1,3 & 2,3 & 1,2,3  \\ \hline
|{\cal X(J)}| & 6 & 12 & 8 & 18 & 14 & 20 & 26  \\ \hline \hline
p=1 & 3^S & 3^S & 3^S & 5 & 5 & 3 & 7 \\ \hline
p=2 & -- & -- & -- & 5 & 5^T & 3 & 5 \\ \hline
p=3 & -- & -- & -- & -- & -- & -- & 7^T \\ \hline \hline
\end{array}$$

$$\begin{array}{||c||c|c|c|c|c|c|c|c|c|c|c|c|c|c|c||}
\hline \hline 
            J & 1 & 2 & 3   & 4 & 1,2 & 1,3 & 1,4 & 2,3 & 2,4 & 3,4 & 1,2,3 & 1,2,4 & 1,3,4 & 2,3,4 & 1,2,3,4 \\ \hline   
|{\cal X(J)}| & 8 & 24 & 32 & 16 & 32 & 40 & 24   & 56  & 40  & 48  & 64    & 48    & 56    & 72 & 80    \\ \hline \hline
p=1 & 3^S & 5^S & 3^S & 3^S & 3 & 5 & 5^S & 3 & 3 & 3 & 7 & 7 & 5 & 5 & 5 \\ \hline
p=2 & -- & -- & -- & -- & 3 & 5 & 5 & 3 & 3 & 3 & 7 & 7^T & 5 & 5 & 7 \\ \hline
p=3 & -- & -- & -- & -- & -- & -- & -- & -- & -- & -- & 5 & 5 & 5 & 5 & 7 \\ \hline \hline
\end{array}$$

Entries marked by an ``S'' indicate \emph{spherical designs}, that is, those with $p=1$ and a constant weight function.  This is obvious when $|J|=1$; for example, for $n=4$, ${\cal I}^4_2$ is the spherical 5-design formed by vectors with minimal non-zero norm in the lattice $D_4$.  Additionally, for $n=4$ we have ${\cal I}^4_1 \cup \frac{1}{2} {\cal I}^4_4$ as a spherical 5-design: it is the set of vectors with minimal non-zero norm in the dual lattice $D_4^*$.  That this is the case follows from the last sentence of the statement of Proposition \ref{L} since for $n=4$ we have $n^2/2^n=1$.

Entries marked by a ``T'' indicate \emph{antipodal tight Euclidean designs}, that is, those of size $N(n,p,t)$.  Indeed, we have  $$N(3,2,5)=14=|{\cal I}^3_1|+|{\cal I}^3_3|,$$ \ $$N(3,3,7)=26=|{\cal I}^3_1|+|{\cal I}^3_2|+|{\cal I}^3_3|,$$ and  $$N(4,2,7)=48=|{\cal I}^4_1|+|{\cal I}^4_2|+|{\cal I}^4_4|.$$  

Our last task is to consider these tight cases.

\emph{Proof of Corollary \ref{FG}.}  Let first $n=3$ and $t=5$.  We then have $${\cal X}=r_1 {\cal I}^3_1 \cup \frac{r_3}{\sqrt{3}} {\cal I}^3_3.$$    
Equation (\ref{M5}) then becomes
$$4r_1^4w_1-\frac{32}{9} r_3^4w_3=0,$$
and this holds exactly when $$w_3=\frac{9}{8} \cdot \frac{r_1^4}{r_3^4} \cdot w_1,$$ proving Corollary \ref{FG}, 1.

Next, let $n=3$ and $t=7$.  Now we have $${\cal X}={\cal I}^3_1 \cup \frac{r_2}{\sqrt{2}} {\cal I}^3_2\cup \frac{r_3}{\sqrt{3}} {\cal I}^3_3.$$  This time we use equations (\ref{DB}), which become    
$$\left\{ 
\begin{array}{l}
4r_1^4w_1-2 r_2^4w_2-\frac{32}{9} r_3^4w_3 =0 \\ \\
 12r_1^6w_1-39r_2^6w_2+\frac{256}{9} r_3^6w_3 =0 \\ \\
4r_1^6w_1-2 r_2^6w_2-\frac{32}{9} r_3^6w_3 =0 
\end{array}
\right.$$
All three equations hold, if and only if, $$r_2=\sqrt{\frac{2r_3^2}{5r_3^2-3r_3^2}} \cdot r_1, \mbox{   } w_2=\frac{1}{10} \cdot \frac{(5r_3^2-3r_1^2)^3}{r_3^6} w_1, \mbox{   and   } w_3=\frac{27}{40} \cdot \frac{r_1^6}{r_3^6} \cdot w_1.$$
Letting $r=r_1$, $\rho=\sqrt{\frac{5r_3^2-3r_1^2}{2}}$, and $w=w_1$ produces Corollary \ref{FG}, 2.

Finally, let $n=4$ and $t=7$.  In this case we have $${\cal X}={\cal I}^4_1 \cup \frac{r_2}{\sqrt{2}} {\cal I}^4_2\cup \frac{r_4}{2} {\cal I}^4_4.$$  The system of equations (\ref{DB}) becomes
$$\left\{ 
\begin{array}{l}
4r_1^4w_1-4 r_4^4w_4 =0 \\ \\
 12r_1^6w_1-36r_2^6w_2+24r_4^6w_4 =0 \\ \\
4r_1^6w_1-4r_4^6w_4 =0 
\end{array}
\right.$$
All three equations hold, if and only if, $r_1=r_4$, $w_1=w_4$, and $w_2=\frac{r_1^6}{r_2^6} \cdot w_1$, as stated in Corollary \ref{FG}, 3.  This completes the proof of Corollary \ref{FG}.

\emph{Remarks.}  Bannai \cite{Ban:Prepb} proved that, up to similarity, the design in Corollary \ref{FG}, 1 is one of exactly three different antipodal tight 5-designs for $p=2$ and $n \geq 3$ (the other two are not fully symmetric).  In Corollary \ref{FG}, 2, two particularly nice choices for the parameters are worth mentioning: with $r=1$ and $\rho=\sqrt{\frac{8}{3}}$, we have
$${\cal X}={\cal I}^3_1 \cup {\cal I}^3_2\cup \frac{1}{2} {\cal I}^3_3$$ with weights $w_2=\frac{1}{10} w$ and $w_3=\frac{8}{5}w$; and with $r=1$ and $\rho=\sqrt{6}$, we have
$${\cal X}={\cal I}^3_1 \cup \frac{1}{2}{\cal I}^3_2\cup {\cal I}^3_3$$ with weights $w_2=\frac{32}{5}w$ and $w_3=\frac{1}{40}w$.  Finally, in Corollary \ref{FG}, 3, we see that with $r=1$ and $\rho=\sqrt{2}$ we have $${\cal X}=({\cal I}^4_1 \cup \frac{1}{2} {\cal I}^4_4) \cup  {\cal I}^4_2,$$
where ${\cal I}^4_2$ is the set of vectors with minimal non-zero norm in the ``checkerboard'' lattice $D_4$, and ${\cal I}^4_1 \cup \frac{1}{2} {\cal I}^4_4$ is the set of vectors with minimal non-zero norm in its dual lattice $D_4^*$.  De la Harpe, Pache, and Venkov \cite{DelPacVen:Prepb} used the root lattice $D_n$ and its dual $D_n^*$ to construct (non-tight) Euclidean $t$-designs for $p=1$ with several higher values of $n$ and $t$.

{\bf Acknowledgments.}  Work on this paper began during a workshop at the University of Geneva, Switzerland, organized by Pierre de la Harpe during May 2004; and it was continued at a second workshop at Kyushu University, in Fukuoka, Japan, organized by Eiichi Bannai and Etsuko Bannai during January, 2005.  I am very grateful for their hospitality and the many helpful discussions.


\begin{thebibliography}{10}





\bibitem{Baj:1992a}
B.~Bajnok.
\newblock Construction of spherical $t$-designs.
\newblock {\em Geom. Dedicata}, 43/2:167--179, 1992.

\bibitem{Baj:2005a}
B.~Bajnok.
\newblock {O}n {E}uclidean designs.
\newblock {\em Adv. Geom.} 6:431--446, 2006.

\bibitem{Ban:1988a}
Ei.~Bannai.
\newblock On extremal finite sets in the sphere and other metric spaces.
\newblock {\em London Math. Soc. Lecture Note Ser.}, 131:13--38, 1988.

\bibitem{BanBan:Prepa}
Ei.~Bannai and Et.~Bannai.
\newblock On Euclidean tight 4-designs.
\newblock Preprint.

\bibitem{Ban:Prepb}
Et.~Bannai.
\newblock On antipodal Euclidean tight $(2e+1)$-designs.
\newblock Preprint.

\bibitem{ConSlo:1999a}
J. H. Conway and N. J. A. Sloane.
\newblock {\em Sphere Packings, Lattices and Groups, 3rd ed}.
\newblock Springer--Verlag, 1999.

\bibitem{DelPac:Prepa}
P. de la Harpe and C. Pache.
\newblock Cubature formulas, geometrical designs, reproducing kernels, and Markov operators.
\newblock Preprint.

\bibitem{DelPacVen:Prepb}
P. de la Harpe, C. Pache, and B. Venkov.
\newblock Construction of spherical cubature formulas using lattices.
\newblock Preprint.

\bibitem{DelGoeSei:1977a}
P.~Delsarte, J.~M. Goethals, and J.~J. Seidel.
\newblock Spherical codes and designs.
\newblock {\em Geom. Dedicata}, 6/3:363--388, 1977.

\bibitem{DelSei:1989a}
P.~Delsarte and J.~J. Seidel.
\newblock Fisher type inequalities for Euclidean $t$-designs.
\newblock {\em Linear Algebra Appl.}, 114--115:213--230, 1989.

\bibitem{Erd:1953a}
A. Erd\'elyi et al.
\newblock Higher transcendental functions, Vol. II. (Bateman Manuscript Project).
\newblock MacGraw--Hill, 1953.

\bibitem{EriZin:2001a}
T. Ericson and V. Zinoviev.
\newblock Codes on Euclidean spheres.
\newblock Elsevier, 2001.

\bibitem{God:1993a}
C.~D. Godsil.
\newblock {\em Algebraic Combinatorics}.
\newblock Chapman and Hall, Inc., 1993.

\bibitem{GoeSei:1979a}
J.~M. Goethals and J.~J. Seidel.
\newblock Spherical designs.
\newblock In D.~K. Ray-Chaudhuri, editor, {\em Relations between combinatorics
  and other parts of mathematics}, 
\newblock {\em Proc. Sympos. Pure Math.}, 34:255--272. Amer. {M}ath. {S}oc., {P}rovidence, {RI}, 1979.

\bibitem{GoeSei:1981a}
J.~M. Goethals and J.~J. Seidel.
\newblock Cubature formulae, polytopes and spherical designs.
\newblock In C.~Davis, B.~Gr{\"{u}}nbaum, and F.~A. Sher, editors, {\em The
  Geometric Vein: The Coxeter Festschrift}, pp. 203--218. Springer-{V}erlag
  New York -- Berlin, 1981.

\bibitem{Hil:1974a}
F. B. Hildebrand.
\newblock {\em Introduction to Numerical Analysis, 2nd ed}.
\newblock McGraw--Hill, Inc., 1974.

\bibitem{NeuSei:1988a}
A. Neumaier and J. J. Seidel.
\newblock Discrete measures for spherical designs, eutactic stars and lattices.
\newblock {\em Nederl. Akad. Wetensch. Indag. Math.}, 50/3:321-334, 1988.

\bibitem{Sei:1990a}
J.~J. Seidel.
\newblock Designs and approximation.
\newblock {\em Contemp. Math.}, 111:179--186, 1990.

\bibitem{Sei:1994a}
J.~J. Seidel.
\newblock Isometric embeddings and geometric designs.
\newblock {\em Discrete Math.}, 136:281--293, 1994.

\bibitem{Sei:1996a}
J.~J. Seidel.
\newblock Spherical designs and tensors.
\newblock In E.~Bannai and A.~Munemasa, editors, {\em Progress in algebraic
  combinatorics}, 
\newblock {\em Adv. Stud. Pure Math.}, 24:309--321. Math. {S}oc. {J}apan, {T}okyo, 1996.

\bibitem{SeyZas:1984a}
P.~D. Seymour and T.~Zaslavsky.
\newblock Averaging sets: A generalization of mean values and spherical designs.
\newblock {\em Adv. Math.}, 52:213--240, 1984.

\bibitem{Str:1971a}
A. H. Stroud.
\newblock {\em Approximate Calculation of Multiple Integrals}.
\newblock Prentice--Hall, Inc., 1971.


\end{thebibliography}
\end{document}